\documentclass[11pt]{article}
\usepackage{latexsym}
\usepackage{amssymb}
\usepackage{amsmath}
\newtheorem{thm}{Theorem}[section]

\newtheorem{lem}[thm]{Lemma}

\newtheorem{prop}[thm]{Proposition}

\newtheorem{definition}[thm]{Definition}

\newtheorem{Bsp.}{Example}[section]

\newtheorem{def.}[theorem]{Definition}
\newtheorem{prop.}[theorem]{Proposition}
\newtheorem{lem.}[theorem]{Lemma}
\newtheorem{cor.}[theorem]{Corollary}
\newtheorem{conj.}[theorem]{Conjecture}
\newtheorem{rem.}{Remarks}[theorem]
\newtheorem{thm.}[theorem]{Theorem}

\newcommand{\h}{\mathcal{H}}
\newcommand{\K}{\mathcal{K}}

\newcommand{\cL}{\mathcal{L}}
\newcommand{\cG}{\mathcal{G}}

\begin{document}


\vspace*{.5cm}
\begin{center}
{\bf   Controlled G-Frames and Their G-Multipliers in Hilbert
spaces}
\bigskip

{\large
{\bf Asghar Rahimi$^a$, Abolhassan Fereydooni$^b$}
\\
{$^a$Department of Mathematics}
\\
 { University of Maragheh}
 \\
    Maragheh, Iran }
\\
{\tt rahimi@maragheh.ac.ir,  } \\
{\large
{$^b$Department of Mathematics}
\\
 {Vali-e-Asr University of Rafsanjan}
 \\
    Rafsanjan, Iran}
\\
{\tt Fereydooniman@yahoo.com  }
\end{center}
\bigskip

\begin{abstract}

 Multipliers have been recently introduced by P. Balazs as operators for Bessel sequences and frames in Hilbert spaces.
 These  are operators that combine (frame-like) analysis, a multiplication
with a fixed sequence ( called the symbol) and synthesis. Weighted
and controlled frames have been introduced to improve the
   numerical efficiency of iterative algorithms for inverting the frame operator
        Also $g$-frames are the most popular  generalization of frames that include almost all of the frame extensions.
     \par In this manuscript the concept of the controlled $g$-frames will be defined and we will show that controlled $g$-frames are equivalent
     to $g$-frames and so the controlled operators $C$ and $C'$ can be used as preconditions in applications.
     Also the multiplier operator for this
      family of operators  will be introduced and some of its properties will be shown.
\end{abstract}

\section{Introduction}
\noindent
\vskip 0.8 true cm
 In \cite{1960 Schatten Norm Ideals of Completely Continious (schatt1)}, R. Schatten provided a detailed study of ideals of compact operators by using their singular decomposition. He investigated the operators of the form $\sum_{k}\lambda_{k}\varphi_{k}\otimes\overline{\psi_{k}}$ where $(\phi_{k})$ and $(\psi_{k})$ are orthonormal families. In \cite{2007 Balazs Basic definition and (1)} the orthonormal families were replaced with Bessel and frame sequences to define Bessel and frame multipliers.
\begin{definition}
Let $\mathcal{H}_{1}$ and $\mathcal{H}_{2}$ be Hilbert spaces, let
$(\psi_{k})\subseteq\mathcal{H}_{1}$ and
$(\phi_{k})\subseteq\mathcal{H}_{2}$ be Bessel sequences. Fix
$m=(m_k)\in l^{\infty}$. The operator ${\bf M}_{m, (
\phi_k),(\psi_k)} : \mathcal{H}_{1} \rightarrow \mathcal{H}_{2}$
defined by
$$ {\bf M}_{m, (\phi_k), ( \psi_k )} (f)  =  \sum \limits_k m_k \langle f,\psi_k\rangle \phi_k $$
is called the \textbf{Bessel multiplier} for the Bessel sequences $(\psi_{k})$ and $(\phi_{k})$. The sequence $m$ is called the symbol of {\bf M}.
\end{definition}
Several basic properties of these operators were investigated in \cite{2007 Balazs Basic definition and (1)}.
 Multipliers are not only interesting from a theoretical point of view, see e.g.
  \cite{2007 Balazs Basic definition and (1),2003 Feichtinger A first,2009 orfler and Torresani Representation of operators in the (doetor09)},
  but they are also used in applications, in particular in the field of audio and acoustic. They have been investigated
 for fusion frames \cite{2008 Arias Bessel fusion},
for generalized frames \cite{2011 A  Rahimi Multipliers of
Genralized frames (5)}, $p$-frames in Banach spaces \cite{2010
Rahimi and P Balazs Multipliers of p-Bessel sequences in (AsPe)}
and for Banach frames \cite{2012 Banach Pair Frames- A. Fereydooni
A. Safapour,2011 Pair Frames- A. Fereydooni A. Safapour}. In
signal processing they are used for Gabor frames
  under the name of Gabor filters \cite{2002 Matz and F Hlawatsch inear Time-Frequency Filters (hlawatgabfilt1)},
  in computational auditory scene analysis they are known by the name of time-frequency masks \cite{2007  Matz D Schafhuber Karlheinz Analysis Optimization and Implementation of (grhahlmasc07)}.
 In real-time implementation of filtering system they approximate time-invariant filters \cite{2005 Balazs WA Deutsch NollSTx Programmer Guide (stx05)}. As a particular way to implement time-variant filters they are used for example for sound morphing \cite{2007 Depalle Kronland-Martinet Time-frequency multipliers for (DepKronTor07)} or psychoacoustical modeling \cite{2009 Balazs Laback Eckel Introducing time-frequency sparsity by  (xxllabmask1)}.

$G$-frames, introduced by W. Sun in \cite{2006 Sun G-frames and
G-Riesz bases (WS)} and improved by the first author \cite{2011 A
Rahimi Frames and Their Generalizations in Hilbert and Banach
Spaces}, are a natural generalization of frames which cover many
other extensions of frames, e.g. bounded quasi-projectors
\cite{2003 Fornasier ecomposition of Hilbert spaces (Fa1),2004
Fornasier Quasi-orthogonal decompsition of structur (Fa2)},
pseudo-frames \cite{2004 Li and H Ogawa Pseudoframes for subspaces
with (Li)}, frame of subspaces or fusion frames \cite{2004 Casazza
Kutyniok Frames of Subspaces (ck)}, outer frames \cite{2004
Aldroubi C CabrelliWavelets on irregular grids (Ald)}, oblique
frames \cite{2004 Christensen Eldar Oblique dual frames and
shift(chel),2003 Eldar  Sampling with arbitrary sampling and
reconstr (eldar)},  and a class of time-frequency localization
operators \cite{2003 orfler Gabor analysis for a class of signals
(Dorf)}. Also it was shown that $g$-frames are equivalent to
stable spaces splitting studied in \cite{1994 Oswald Multilevel
Finite Element Approximation (osw)}. All of these concepts are
proved to be useful in many applications. Multipliers for
$g$-frames introduced in \cite{2011 A  Rahimi Multipliers of
Genralized frames (5)} and some of its properties investigated.

Weighted and controlled frames have been introduced recently to improve the numerical efficiency of iterative
algorithms for inverting the frame operator on abstract
 Hilbert spaces \cite{1010 Balazs J-P Antoine Wighted and Controlled Frames (peter weight)},
 however they are used earlier in \cite{2005 Bogdanova  Vandergheynst Antoine Stereographic wavelet frames (bogd)}
  for spherical wavelets. In this manuscript the concept of controlled  $g$- frame will be defined and we will show that
  any controlled $g$-frame is equivalent a $g$-frame and the role of controller operators are like the role of preconditions
  matrices or operators in linear algebra. Furthermore the multiplier
   operator for these family will be investigated.
\par The paper is organized as follows. In Section 2 we fix the
notation in this paper, summarize known and prove some new results
needed for the rest of the paper. In Section 3 we will define the
concept of controlled $g$-frames and we will show that a
controlled $g$-frame is equivalent to a $g$-frame and so the
controlling operators can be used as precondition matrices in the
problems related to applications.
 In section 4 we will define multipliers of controlled $g$-frame operators and we will prove some of its properties.

\section{Preliminaries}

Now we state some notations and theorems which are used in the present paper. Through this paper, $\h$ and $\K$ are Hilbert spaces and $\{\,\h_{i}:{i\in I}\}$ is a sequence of Hilbert spaces, where $I$ is a subset of $\mathbb{Z}.$ $\cL(\h,\K_{}) $ and $\cL(\h)$ is the collection of all bounded linear operators from $\h$ into $\K $ and $\h$ respectively.

A bounded operator $T$ is called \emph{positive} (respectively \emph{non-negative}), if $\langle Tf,f\rangle>0$ for all $f\neq 0$ (respectively $\langle Tf,f\rangle\geq0$ for all $f$). Every non-negative operator is clearly self-adjoint.
For $T_{1},T_{2} \in \cL(\h)$, we write $ T_{1} \leqslant T_{2}  $ whenever
$$ \langle T_{1}(f) , f \rangle \leqslant\langle T_{2}(f) , f \rangle , \quad \forall f \in \h .$$

If $U\in \cL(\h)$ is non-negative, then there exists a unique non-negative operator $V$ such that $V^2=U$. Furthermore $V$ commutes with every operator that commute with $U$. This will be denoted by $V=U^{\frac{1}{2}}$. Let $\cG\cL(\h)$ be the set of all bounded operators with a bounded inverse and $\cG\cL^+(\h)$ be the set of positive operators in $\cG\cL(\h)$. For $U \in \cL(\h)$, $U \in \cG\cL^+(\h)$ if and only if there exists $0<m \leqslant M<\infty$ such that
$$ m \leqslant U \leqslant M .$$
For $U^{-1}$ we have
$$  M^{-1} \leqslant U^{-1} \leqslant m^{-1}   .$$
The following theorem can be found in \cite{1982 Heuser Functional
analysis}.
\begin{thm}\label{t:order of positive operators}
Let $T_{1},T_{2}, T_{3} \in \cL(\h)$ and  $ T_{1} \leq T_{2} $.
Suppose $T_{3} > 0$ commutes with $T_{1}$ and $T_{2}$ then
$$ T_{1} T_{3}  \leq T_{2} T_{3} .$$
\end{thm}
\par Recall that if $T$ is a compact operator on a separable
Hilbert space $\mathcal{H}$, then in \cite{1960 Schatten Norm
Ideals of Completely Continious (schatt1)} it is proved that there
exist orthonormal sets $\{e_{n}\}$ and $\{\sigma_{n}\}$ in
$\mathcal{H}$ such that
$$Tx=\sum_{n}\lambda_{n}\langle x, e_{n}\rangle\sigma_{n},  $$
for $x\in \mathcal{H}$, where $\lambda_{n}$ is the $n$-th singular
value of $T$. Given $0<p<\infty$, the \textbf{\textbf{Schatten
$p$-class }}of $\mathcal{H}$ \cite{1960 Schatten Norm Ideals of
Completely Continious (schatt1)}, denoted $\mathcal{S}_{p}$, is
the space of all compact operators $T$ on $\mathcal{H}$ with the
singular value sequence $\{\lambda_{n}\}$ belonging to $\ell^{p}$.
It was shown that \cite{1990 Zhuo Operator Theory in Function
Spaces(Zhuo) } , $\mathcal{S}_{p}$ is a Banach space with the norm
\begin{equation}
\label{normp}\|T\|_{p}=[\sum_{n}|\lambda_{n}|^{p}]^{\frac{1}{p}}.
\end{equation}
$\mathcal{S}_{1}$ is called the\textit{ trace class} of
$\mathcal{H}$ and $\mathcal{S}_{2}$ is called the
\textit{Hilbert-Schmidt class.} $T\in \mathcal{S}_p$ if and only
if $T^{p}\in \mathcal{S}_1$. Moreover
$\|T\|^{p}_{p}=\|T^{p}\|_{1}$. Also, $T\in \mathcal{S}_{p}$ if and
only if $|T|^p=(T^*T)^\frac{p}{2}\in \mathcal{S}_{1}$ if and only
if $T^{*}T\in \mathcal{S}_{\frac{p}{2}}$. Moreover,
$\|T\|^{p}_{p}=\|T^{*}\|^{p}_{p}=\||T|\|^{p}_{p}=\||T|^{p}\|_{1}=\|T^{*}T\|_{\frac{p}{2}}.$

\vskip 0.8 true cm

\subsection{\bf $G$-Frames} \noindent \vskip 0.8 true cm
For any sequence $\{\,\h_{i}:{i\in I}\},$ we can assume that there exits a Hilbert space $\K$ such that for all $i\in I, \h_{i}\subseteq \K$
 ( for example $\K={(\bigoplus_{i\in I}\h_{i})}_{\ell^2}$ ).
\vskip 0.4 true cm
\begin{def.}
A sequence $\Lambda=\{\,\Lambda_{i}\in \cL(\h, \h_{i}): i\in I\,\}$ is called generalized  frame, or simply a \textbf{$g$-frame}, for $\h$ with respect to $\{\h_{i}: i\in I\}$  if there exist constants $A>0$ and $B<\infty$ such that
\begin{equation}\label{deframe}
A\|f\|^{2}\leqslant \sum_{i\in I}\| \Lambda_{i}f\|^{2} \leqslant B\|f\|^{2},
\quad \forall f\in \h.
\end{equation}
The numbers $A$ and $B$ are called $g$-frame bounds.
\end{def.}
\vskip 0.4 true cm
 $\Lambda=\{\Lambda_{i}: i\in I\}$ is called  \emph{tight} $g$-frame if $A = B$ and \emph{Parseval} $g$-frame if $A=B=1$. If the second inequality in (\ref{deframe}) holds, the sequence is called $g$-Bessel sequence.
\vskip 0.4 true cm
$\Lambda=\{\,\Lambda_{i}\in \cL(\h, \h_{i}): i\in I\,\}$ is called a \textbf{$g$-frame sequence}, if it is a $g$-frame for $\overline{\text{span}}\{\Lambda_{i}^{*}(\h_{i})\}_{i\in I}$.
 \vskip 0.4 true cm
 It is easy to see that, if $\{f_i\}_{i\in I}$ is a frame for $\h$ with bounds $A$ and $B$, then by putting $\h_i=\mathbb{C}$ and $\Lambda_{i}(\cdot)=\langle \cdot , f_i\rangle$, the family $\{\Lambda_{i}: i\in I\}$ is a $g$-frame for $\h$ with bounds $A$
and $B$.
 Let
\begin{equation}\label{e:20 }
 \Big(\bigoplus_{i\in I}\h_{i}\Big)_{\ell_{2}}=\Big\{\,\{f_{i}\}_{i\in I}\mid f_{i}\in \h_{i},\,\forall i\in I\ \ \text{and}\ \  \sum_{i\in I}\|
 f_{i}\|^{2}<+\infty\Big\}.
 \end{equation}
\begin{prop.}\label{gbessel}\cite{2009 Najati and A Rahimi Generalized Frames in Hilbert spaces (NaRa)}
$\Lambda=\{\,\Lambda_{i}\in \cL(\h, \h_{i}): i\in I\,\}$ is a $g$-Bessel sequence for $\h$ with bound $B$, if and only if the operator
\[T_{ \Lambda }:\Big(\bigoplus_{i\in I}\h_{i}\Big)_{\ell_{2}}\longrightarrow \h\]
defined by
\[T_{ \Lambda}(\{f_{i}\}_{i\in I})=\sum_{i\in I}\Lambda_{i}^{*}(f_{i})\]
is a well-defined and bounded operator with $\|T_{\mathbf{\Lambda}}\|\leqslant \sqrt{B}$. Furthermore
\[T^{*}_{ \Lambda}:\h\longrightarrow\Big(\bigoplus\h_{i}\Big)_{\ell_{2}}\]
\[T^{*}_{ \Lambda}(f)=\{\Lambda_{i}f\}_{i\in I}.\]
\end{prop.}
\vskip 0.4 true cm
If $\Lambda=\{\,\Lambda_{i}\in \cL(\h, \h_{i}): i\in I\,\}$ is a $g$-frame, the operators $T_{ \Lambda}$ and $T^{*}_{ \Lambda}$ in Proposition  \ref{gbessel} are called {\bf synthesis operator} and {\bf analysis operator} of $ \Lambda =\{\,\Lambda_{i}\in \cL(\h, \h_{i}): i\in I\,\},$
respectively.
\vskip 0.4 true cm
\begin{prop.}\label{Tonto}\cite{2008 Najati Rahimi G-frames }
$\Lambda=\{\,\Lambda_{i}\in \cL(\h, \h_{i}): i\in I\,\}$ is a
$g$-frame for $\h$ if and only if  the synthesis operator
$T_{\Lambda}$ is well-defined, bounded and onto.
\end{prop.}
The ordinary version of the next theorem which is proved in \cite{2012 Banach Pair Frames- A. Fereydooni A. Safapour}, can be extended easily to the general case.

\begin{thm}\label{t:g Bessel operator welldefined}
\cite{2012 Banach Pair Frames- A. Fereydooni A. Safapour}
$\Lambda=\{\,\Lambda_{i}\in \cL(\h, \h_{i}): i\in I\,\}$ is a $g$-Bessel sequence for $\h$ if and only if the operator
\begin{equation}\label{e:Bessel operator}
S_{\Lambda}: \h\longrightarrow \h,\quad S_{\Lambda}=\sum_{i\in I}\Lambda^{\ast}_{i}\Lambda_{i}f ,
\end{equation}
is a welldefined operator. In this case $S_{\Lambda}$ is bounded.
\end{thm}
\begin{thm}\label{t:g frame operator welldefined}
\cite{2011 Pair Frames- A. Fereydooni A. Safapour}
$\Lambda=\{\,\Lambda_{i}\in \cL(\h, \h_{i}): i\in I\,\}$ is a $g$-frame for $\h$ if and only if the operator
$$ S_{\Lambda}: \h\longrightarrow \h,\quad S_{\Lambda}=\sum_{i\in I}\Lambda^{\ast}_{i}\Lambda_{i}f ,$$
is a welldefined invertible operator. In this case $S_{\Lambda}$ is bounded.
\end{thm}

$ S_{\Lambda}$ is called the \textbf{$g$-frame operator }of $\Lambda=\{\,\Lambda_{i}: i\in I\,\}$ and it is known \cite{2009 Najati and A Rahimi Generalized Frames in Hilbert spaces (NaRa)} that $S_{\Lambda}$ is a positive and
\[AI\leqslant S_{\Lambda}\leqslant BI,\]
where A and B are the frame bounds. Every $f\in\h$ has an
expansion $f=\sum_{i}\Lambda_{i}^{*}\Lambda_{i}S^{-1}_{\Lambda}f$.
One of the most important advantages  of $g$-frames is  a
resolution of identity $\sum_{i}\Lambda_{i}^{*}\Lambda_{i}S^{-1}_{
\Lambda}=I$. \vskip 0.4 true cm \vskip 0.8 true cm
\subsection{\bf Multipliers of $g$-frames}
\noindent
\vskip 0.8 true cm
 The concept of multipliers for $g$-Bessel sequences introduced by the first author in \cite{2011 A  Rahimi Multipliers of Genralized frames (5)} and some of their properties will be shown.

\begin{def.}
Let $\Lambda=\{\,\Lambda_{i}\in \cL(\h, \h_{i}): i\in I\,\}$ and $\Theta=\{\,\Theta_{i}\in \cL(\h, \h_{i}): i\in I\,\}$ be $g$-Bessel sequences. 
If for $m \subset \mathbb{C}$, the operator
\[ \mathbf{M}=\mathbf{M}_{m,\Lambda,\Theta}:\h\rightarrow\h\]
\begin{equation}\label{defmul}
\mathbf{M}(f)=\sum_{i}m_{i}\Lambda^{*}_{i}\Theta_{i}f
\end{equation}
is welldefined, then $\mathbf{M}$ is called the \textbf{$g$-multiplier} of $\Lambda$,$\Theta$ and $m$.
\end{def.}
\vskip 0.4 true cm
 If $m=(m_i)=(1, 1, 1,...)$ and $\mathbf{M}=I$, $(\Lambda,\Theta)$ is called a \emph{pair dual} (i.e. $I=\sum_{i\in I}\Lambda^*_i\Theta_i$ ).
 \vskip0.4 true cm
 Let $\{\lambda_{i}\}$ and $\{\varphi_{i}\}$ be  Bessel sequences and $m\in\ell^{\infty}$, consider the
 corresponding $g$-Bessel sequences $ \Lambda_{i}\cdot=\langle\cdot,\lambda_{i}\rangle $
 and $\Theta_{i}\cdot=\langle\cdot,\varphi_{i}\rangle $. For any $f\in\h$ we have:
\begin{equation*}
\begin{aligned}
\mathbf{M}_{m,\Lambda,\Theta}(f) ={\bf M}_{m, ( \lambda_k ),
(\phi_k) } (f)= \sum_i m_i \langle f ,\varphi_i \rangle \lambda_i.
\end{aligned}
\end{equation*}

It is easy to show that the adjoint of
$\mathbf{M}_{m,\Lambda,\Theta}$ is
$\mathbf{M}_{\overline{m},\Theta,\Lambda}$.

\vskip 0.4 true cm
\begin{lem.}\cite{2011 A  Rahimi Multipliers of Genralized frames (5)}
If  $\Theta=\{\,\Theta_{i}\in \cL(\h, \h_{i}): i\in I\,\}$ is a $g$-Bessel sequence with bound $B_{\Theta}$ and $m=(m_i)\in \ell^{\infty}$, then $\{m_{i}\Theta_{i}\}_{i\in I}$ is a $g$-Bessel sequence with bound $\|m\|_{\infty}B_{\Theta}$.
\end{lem.}
 Like weighted frames \cite{1010 Balazs J-P Antoine Wighted and Controlled Frames (peter weight)}, $\{m_{i}\Theta_{i}\}_{i\in I}$ can be called \emph{ weighted $g$- frame ($g$-Bessel)}.  By using the synthesis  and the analysis operators of $\Lambda$ and $ m \Theta$, respectively,  we can write

\[\mathbf{M}_{m,\Lambda,\Theta}f=\sum_{i}m_{i}\Lambda_{i}^{*}\Theta_{i}f=\sum_{i}\Lambda_{i}^{*}(m_{i}\Theta_{i})f =T_{\Lambda}\{m_i \Theta_i f\}=T_{\Lambda}T^{*}_{m\Theta}f.\]
So
\begin{equation}\label{defmul2}
\mathbf{M}_{m,\Lambda,\Theta}=T_{\Lambda}T^{*}_{m \Theta}.
\end{equation}
If we define the diagonal operator
\[ D_m:\Big(\bigoplus\h_{i}\Big)_{\ell_{2}}\rightarrow\Big(\bigoplus\h_{i}\Big)_{\ell_{2}},\]
\begin{equation}\label{defmul21}
D_m((\xi_{i}))=(m_i \xi_i)_{i\in I}
\end{equation}
then
\begin{equation}\label{defmul3}
\mathbf{M}_{m,\Lambda,\Theta}=T_{\Lambda}D_m T^{*}_{\Theta}.
\end{equation}
The  notations in (\ref{defmul2}), (\ref{defmul21}) and
(\ref{defmul3}) were used for proving the following propositions
in \cite{2011 A  Rahimi Multipliers of Genralized frames (5)} .
\begin{prop.}\cite{2011 A  Rahimi Multipliers of Genralized frames (5)}
Let $m\in \ell^{\infty}$, $\Lambda=\{\,\Lambda_{i}\in \cL(\h, \h_{i}): i\in I\,\}$ be a $g$-Riesz base
 and $\Theta=\{\,\Theta_{i}\in \cL(\h, \h_{i}): i\in I\,\}$ be a $g$-Bessel sequence.
 The map $m\rightarrow \mathbf{M}_{m,\Lambda,\Theta}$ is injective.
\end{prop.}
\begin{prop.}\cite{2011 A  Rahimi Multipliers of Genralized frames (5)}
Let $\Lambda=\{\,\Lambda_{i}\in \cL(\h, \h_{i}): i\in I\,\}$ and $\Theta=\{\,\Theta_{i}\in \cL(\h, \h_{i}): i\in I\,\}$
 be $g$-Bessel sequences for $\h$. If $m=(m_i)\in c_0$ and $(rank \Theta_{i})\in \ell^{\infty}$,
 then $\mathbf{M}_{m,\Lambda,\Theta}$ is compact.
 \end{prop.}

\begin{prop.}\cite{2011 A  Rahimi Multipliers of Genralized frames (5)}
Let $\Lambda=\{\,\Lambda_{i}\in \cL(\h, \h_{i}): i\in I\,\}$ and $\Theta=\{\,\Theta_{i}\in \cL(\h, \h_{i}): i\in I\,\}$ be $g$-Bessel sequences for $\h$. If $m=(m_i)\in \ell^{p}$ and $(dim\h_i)_{i\in I}\in\ell^{\infty}$, then $\mathbf{M}_{m,\Lambda,\Theta}$ is a Schatten $p$-class operator.
\end{prop.}
\vskip 0.4 true cm

\begin{cor.}\cite{2011 A  Rahimi Multipliers of Genralized frames (5)}
Let $\Lambda=\{\,\Lambda_{i}\in \cL(\h, \h_{i}): i\in I\,\}$ and $\Theta=\{\,\Theta_{i}\in \cL(\h, \h_{i}): i\in I\,\}$ be $g$-Bessel sequences for $\h$.
\begin{enumerate}
 \item If $m=(m_i)\in \ell^{1}$ and $(dim\h_i)_{i\in I}\in\ell^{\infty}$, then $\mathbf{M}_{m,\Lambda,\Theta}$ is a trace- class operator.
\item If $m=(m_i)\in \ell^{2}$ and $(dim\h_i)_{i\in I}\in\ell^{\infty}$, then $\mathbf{M}_{m,\Lambda,\Theta}$ is a Hilbert-Schmit operator.
\end{enumerate}
\end{cor.}
\vskip 0.4 true cm
\vskip 0.8 true cm
\section{\bf {\bf \em{\bf Controlled  $g$-frames  }}}
\noindent
\vskip 0.8 true cm
Weighted and controlled frames have been introduced recently to improve the numerical efficiency of iterative algorithms for inverting the frame operator. In \cite{1010 Balazs J-P Antoine Wighted and Controlled Frames (peter weight)}, it was shown that the controlled frames are equivalent to standard frames and it was used in the sense of preconditioning.

 In this section, the concepts of controlled frames and controlled Bessel sequences will be extended to $g$-frames
 and we will show that controlled $g$-frames are equivalent $g$-frames.

\begin{def.}\label{c:C,C'-controlled Bessel sequence}
Let $C,C'\in \cG\cL^+(\h)$. The family $\Lambda=\{\,\Lambda_{i}\in
\cL(\h, \h_{i}): i\in I\,\}$ will be called a $(C,C')$-controlled
$g$-frame for $\h$, if $\Lambda$ is a $g$-Bessel sequence  and
there exists constants $A>0$ and $B<\infty$ such that
\begin{equation}\label{e:(C,C')-controlled g-Bessel}
A\|f\|^{2}\leqslant\sum_{i\in I} \langle \Lambda_{i}Cf,
\Lambda_{i}C^{'} f \rangle \leqslant B\|f\|^{2}, \quad \forall
f\in \h.
\end{equation}

A and B will be called \emph{controlled frame bounds}. If $C'=I$,
we call $\Lambda=\{\,\Lambda_{i}\}$ a \textbf{$C$-controlled
$g$-frame} for $\h$ with bounds $A$ and $B$. If the second part of
the above inequality holds, it will be called
\textbf{$(C,C')$-controlled $g$-Bessel sequence} with bound $B$.
 \end{def.}
 The proof of the following lemmas is straightforward.
 \begin{lem}
Let $C\in\cG\cL^{+}(\h)$. The $g$-Bessel sequence
$\Lambda=\{\,\Lambda_{i}\in \cL(\h, \h_{i}): i\in I\,\}$ is
$(C,C)$-controlled Bessel sequence (or  $(C,C)$-controlled
$g$-frame) if and only if there exists constant $B<\infty$ (and
$A>0$ )such that
$$ \, \sum_{i\in I}  \| \Lambda_{i}Cf \|^2 \leqslant B\|f\|^{2}, \quad \forall f\in \h$$
$$ (\, or \,  A\|f\|^{2}\leqslant \sum_{i\in I}  \| \Lambda_{i}Cf \|^2 \leqslant B\|f\|^{2}, \quad \forall f\in \h).$$
\end{lem}
We call the $(C,C)$-controlled Bessel sequence  and
$(C,C)$-controlled $g$-frame, $C^2$-controlled Bessel sequence and
$C^2$-controlled $g$-frame with bounds $A,B$.

\begin{lem}
For $C,C'\in \cG\cL^+(\h)$, the family $\Lambda=\{\,\Lambda_{i}\in
\cL(\h, \h_{i}): i\in I\,\}$ is a $(C,C')$-controlled Bessel
sequence for $\h$ if and only  if the operator
$$ L_{C \Lambda C'}: \h \rightarrow \h, \quad L_{C \Lambda C'} f:=\sum_{i\in I}C'\Lambda_{i}^{*}\Lambda_{i}Cf, $$
 is well defined and there exists constant $B<\infty$ such that
$$\sum_{i\in I} \langle \Lambda_{i}Cf, \Lambda_{i}C^{'} f \rangle
\leqslant B\|f\|^{2}, \quad \forall f\in \h.$$
 \end{lem}
 The operator $$ L_{C \Lambda C'}: \h \rightarrow \h, \quad L_{C \Lambda C'} f:=\sum_{i\in I}C'\Lambda_{i}^{*}\Lambda_{i}Cf,
 $$ is called the $(C,C')$-controlled Bessel
sequence operator, also $ L_{C \Lambda C'} = C S_{\Lambda} C'$. It
follows from the definition that for a $g$-frame, this operator is
positive and invertible and
$$ AI\leq L_{C \Lambda C'}\leq BI.$$
Also, if $C$ and $C'$ commute with each other, then $C' ,C'^{-1},
C ,C^{-1}$ commute with $ L_{C \Lambda C'}, CS_{\Lambda},
S_{\Lambda}C'$ .
\par The following proposition shows that any $g$-frame is a
controlled $g$-frame and versa. This is the most important
advantage of weighted and controlled $g$-frame in the sense of
precondition.
 \begin{prop}
Let $C \in \cG\cL^{+}(\h)$. The family
$\Lambda=\{\,\Lambda_{i}\in \cL(\h, \h_{i}): i\in I\,\}$ is a
$g$-frame if and only if $\Lambda$ is a $C^2$-controlled
$g$-frame.
 \end{prop}
{\bf Proof.} Suppose that $\Lambda$ is a $C^2$-controlled
$g$-frame with bounds $A,B$. Then
$$ \,A\|f\|^{2} \leqslant \sum_{i\in I}  \| \Lambda_{i}C  f \|^2 \leqslant B\|f\|^{2}, \quad \forall f\in \h ,$$
 For $f \in \h$
\[ \begin{split}
A\|f\|^{2} = A\| CC^{-1} f\|^{2} \leqslant A\| C \|^{2} \| C^{-1} f\|^{2} \leqslant \| C \|^{2} \sum_{i\in I}  \| \Lambda_{i}C  C^{-1} f \|^2 \\
= \| C \|^{2} \sum_{i\in I}  \| \Lambda_{i} f \|^2.
\end{split} \]
Hence
$$ A \| C \|^{-2} \|f\|^{2} \leqslant \sum_{i\in I}  \| \Lambda_{i} f \|^2, \quad \forall f \in \h .$$
On the other hand for every $f \in \h$,
\[ \begin{split}
 \sum_{i\in I}  \| \Lambda_{i} f \|^2=\sum_{i\in I}  \| \Lambda_{i} CC^{-1} f \|^2 \leqslant  B\| C^{-1}f\|^{2}\leqslant  B\| C^{-1}\|^{2} \|f\|^{2} .
\end{split} \]
These inequalities yields that $\Lambda$ is a $g$-frame with
bounds $ A \| C \|^{-2}  , B\| C^{-1}\|^{2} $. For the converse
assume that $\Lambda$ is  $g$-frame with bounds $ A' , B' $. Then
for all $f \in \h $,
$$  A' \|f\|^{2} \leqslant \sum_{i\in I}  \| \Lambda_{i} f \|^2 \leqslant  B' \|f\|^{2}.$$
So for  $f \in \h $,
$$  \sum_{i\in I}  \| \Lambda_{i} Cf \|^2 \leqslant  B' \|Cf\|^{2} \leqslant  B'  \|C \|^{2} \|f\|^{2}.$$
For lower bound, the  $g$-frameness of $\Lambda$ shows that for
any if  $f \in \h $,
\[ \begin{split} A' \|f\|^{2} =  A '  \| C^{-1} C f\|^{2} \leqslant A'  \| C^{-1} \|^{2} \| C f\|^{2}
\leqslant   \| C^{-1} \|^{2}\sum_{i\in I}  \| \Lambda_{i} C f \|^2
.\end{split} \] Therefore  $\Lambda$ is a $C^2$-controlled
$g$-frame with bounds $A'   \| C^{-1} \|^{-2} ,B' \|C \|^{2}$.
$ \Box $  \vskip 0.4 true cm
\begin{prop}
 Assume that $\Lambda=\{\,\Lambda_{i}: i\in I\,\}$ is a
$g$-frame and $C,C' \in  \cG\cL^+(\h)$, which commute with each
other and commute with $S_{\Lambda}$.  Then
$\Lambda=\{\,\Lambda_{i}: i\in I\,\}$ is a $(C,C')$-controlled
$g$-frame.
\end{prop}
{\bf Proof.} Let $\Lambda$ be g-frame with bounds $A,B$ and  $m ,
m'>0, M, M'<\infty$ so that
$$
m \leqslant C \leqslant M, \quad m' \leqslant C' \leqslant M' . $$
Then
$$ m A \leqslant C S_{\Lambda} \leqslant M B ,$$
because $C$ commute with $S_{\Lambda}$. Again $C'$ commutes with
$C S_{\Lambda}$ and then
$$ mm'A \leqslant L_{C \Lambda C'} \leqslant MM'B .$$

$ \Box $  \vskip 0.4 true cm

\section{\bf {\bf \em{\bf   Multipliers of Controlled  $g$-frames}}}
\noindent \vskip 0.8 true cm
Extending the concept of multipliers of frames, in this section,
we will define controlled $g$-frame's multiplier for
$C$-controlled $g$-frames in Hilbert spaces. The definition of
general case $( C, C')$-controlled $g$-frames  goes smooth.

\begin{lem.} \label{l:C2 and C'*2-controlled g-Bessel sequences to  multiplier operator}
Let $C,C'\in \cG\cL^+(\h)$ and  $\Lambda=\{\,\Lambda_{i}\in
\cL(\h, \h_{i}): i\in I\,\}$, $\Theta=\{\,\Theta_{i}\in \cL(\h,
\h_{i}): i\in I\,\}$ be ${C'}^2$ and $C^2$-controlled $g$-Bessel
sequences for $\h$, respectively. Let $m\in\ell^{\infty}$. The
operator
 $$\textbf{M}_{m C \Theta  \Lambda C'}:\h\rightarrow\h$$
 defined by
  $$ \textbf{M}_{m C \Theta  \Lambda C'} f:=\sum_{i\in I}m_i C \Theta_i^*  \Lambda_i C' f $$
is a well-defined bounded operator.
\end{lem.}
{\bf Proof.} Assume  $\Lambda=\{\,\Lambda_{i}\in \cL(\h, \h_{i}):
i\in I\,\}$, $\Theta=\{\,\Theta_{i}\in \cL(\h, \h_{i}): i\in
I\,\}$ be  ${C'}^2$ and $C^2$-controlled $g$-Bessel sequences for
$\h$ with bounds $B,B'$, respectively. For any $f,g\in\h$ and
finite subset $J \subset I$,
\begin{align*}
 \| \sum_{i \in J}m_i C \Theta_i^*  \Lambda_i C' f\|= & \sup_{g \in \h, \|g\|=1} \| \sum_{i \in J}m_i  \langle \Lambda_i C' f,\Theta_i C^* g \rangle\|  \\
& \leqslant \sup_{g \in \h, \|g\|=1}\sum_{i \in J} |m_i | \|   \Lambda_i  C' f\|  \|  \Theta_i C^* g \| \\
& \leqslant \sup_{g \in \h, \|g\|=1}\|m\|_{\infty} \big(\sum_{i \in I} \| \Theta_i C^* g \|^2\big)^{\frac{1}{2}} \big(\sum_{i \in J}\|   \Lambda_i C' f\|^2  \\
&  \leqslant \|m\|_{\infty}\sqrt{BB'}\|f\|
\end{align*}
This shows that $\textbf{M}_{m C \Theta  \Lambda C'}$ is
well-defined and $$\|\textbf{M}_{m C \Theta  \Lambda C'}\|\leq
\|m\|_{\infty}\sqrt{BB'}.$$ $ \Box $ \vskip 0.4 true cm
Above Lemma ia a motivation to define the following definition.

\begin{def.}\label{c:}
Let $C,C'\in \cG\cL^+(\h)$ and  $\Lambda=\{\,\Lambda_{i}\in
\cL(\h, \h_{i}): i\in I\,\}$, $\Theta=\{\,\Theta_{i}\in \cL(\h,
\h_{i}): i\in I\,\}$ be ${C'}^2$ and $C^2$-controlled $g$-Bessel
sequences for $\h$, respectively. Let $m\in\ell^{\infty}$. The
operator
$$ \textbf{M}_{m C \Theta  \Lambda C'}:\h\rightarrow\h $$
 defined by
 \begin{equation}\label{e:(C,C')-controlled multiplier operator}
  \textbf{M}_{m C \Theta  \Lambda C'} f:=\sum_{i\in I}m_i C \Theta_i^*  \Lambda_i C' f,
\end{equation}
 is called the $(C,C')$-controlled multiplier operator with symbol $m$.
\end{def.}

By using representations (\ref{defmul2}) and (\ref{defmul3}), we have
$$\textbf{M}_{m C \Theta  \Lambda C'}=C \textbf{M}_{m  \Theta  \Lambda } C'=CT_ \Theta  D_m T^*_\Lambda C'.  $$
The proof of Proposition 4.7. of \cite{2011 A  Rahimi Multipliers of Genralized frames (5)} shows that if $m=(m_i)\in \ell^{p}$ and $(dim\h_i)_{i\in I}\in\ell^{\infty}$, then the diagonal operator $D_m$ is a Schatten $p$-class operator. Since $\mathcal{S}_p$ is a $*$-ideal of $\mathcal{L}(\h)$ so we have:
\begin{thm.}
Let $\Lambda=\{\,\Lambda_{i}\in \cL(\h, \h_{i}): i\in I\,\}$ and $\Theta=\{\,\Theta_{i}\in \cL(\h, \h_{i}): i\in I\,\}$  be controlled $g$-Bessel sequences for $\h$. If $m=(m_i)\in \ell^{p}$ and $(dim\h_i)_{i\in I}\in\ell^{\infty}$, then $\textbf{M}_{m C \Theta  \Lambda C'} $ is a Schatten $p$-class operator.
\end{thm.}
And
\begin{cor.}
Let $\Lambda=\{\,\Lambda_{i}\in \cL(\h, \h_{i}): i\in I\,\}$ and $\Theta=\{\,\Theta_{i}\in \cL(\h, \h_{i}): i\in I\,\}$  be controlled $g$-Bessel sequences for $\h$.
\begin{enumerate}
 \item If $m=(m_i)\in \ell^{1}$ and $(dim\h_i)_{i\in I}\in\ell^{\infty}$, then $\textbf{M}_{m C \Theta  \Lambda C'}$ is a trace- class operator.
 \item If $m=(m_i)\in \ell^{2}$ and $(dim\h_i)_{i\in I}\in\ell^{\infty}$, then $\textbf{M}_{m C \Theta  \Lambda C'}$ is a Hilbert-Schmit operator.
\end{enumerate}
\end{cor.}
\vskip 0.4 true cm
\vskip 0.8 true cm

\textbf{Acknowledgment}: Some of the results in this paper were
obtained during the first author visited the Acoustics Research
Institute, Austrian Academy of Sciences, Austria, he thanks this
institute for their hospitality.


\newpage

\end{document}